\documentclass[11pt]{amsart}
\usepackage{amssymb,amsthm}

\newtheorem{mr}{Theorem}

\newtheorem{theorem}{Theorem}
\newtheorem{corollary}[theorem]{Corollary}
\newtheorem{proposition}[theorem]{Proposition}
\newtheorem{lemma}[theorem]{Lemma}                   

\theoremstyle{remark}
\newtheorem*{example}{Example}

\newtheorem*{remark}{Remark}
\newtheorem*{question}{Question}
\newtheorem*{ack}{Acknowledgement}

\def\tr{\operatorname{tr}}
\def\R{\mathbb{R}}
\def\C{\mathbb{C}}

\def\N{\mathbb{N}}
\def\md{M(d)}
\def\eps{\varepsilon}
\def\gldc{GL(d)}
\def\pgldc{PGL(d)}
\def\sldc{SL(d)}
\def\OO{\mathcal{O}}
\def\RR{\mathcal{R}}
\def\SS{\mathcal{S}}
\def\norma{\|\mathord{\cdot}\|}   %spacing corresponding to ordinary symbol instead of binary operator

\begin{document}

%%%%%%%%%%%%%%%%%%%%%%%%%%%%%%%%%%%%%%%%%%%%%%%
\title{Inequalities for numerical invariants of sets of matrices}
\author{Jairo Bochi}
\date{June 20, 2002}
\thanks{Suported by CNPQ--Profix.}
%\subjclass[2000]{15A45, 14L24} %16R30  % doesn't work at ArXiv

\address{
IMPA \\
Estr. D. Castorina 110 \\
22460-320 Rio de Janeiro -- Brazil.}
\email{bochi@impa.br}

\begin{abstract}
We prove three inequalities relating some invariants of sets of matrices,
such as the joint spectral radius.
One of the inequalities, in which proof we use geometric invariant theory,
has the generalized spectral radius theorem of Berger and Wang as an immediate corollary.
\end{abstract}

\maketitle

%%%%%%%%%%%%%%%%%%%%%%%%%%%%%%%%%%%%%%%%%%%%%%%%%%%%%%%%

\section{Introduction}

Let $\md$ be the space of $d\times d$ complex matrices.
If $A\in \md$, we indicate by $\rho(A)$
the spectral radius of $A$, that is, the maximum absolute value 
of an eigenvalue of $A$.
Given a norm $\norma$ in $\C^d$,
we endow the space $\md$ with the operator norm
\mbox{$\|A\| = \sup \, \{\|Av\|; \; \|v\|=1\}$.}

\smallskip

For every $A \in \md$ and every norm $\norma$ in $\C^d$, we have
$\rho(A) \le \|A\|$.
On the other hand, there is also a \emph{lower} bound for
$\rho(A)$ in terms of norms:
\begin{equation} \label{e.1st ineq}
\| A^d \| \le C \rho(A) \| A \|^{d-1},
\qquad \text{where $C = 2^d -1$.}
\end{equation}
In particular, if $\rho(A) \ll \|A\|$ then $\|A^d\| \ll \|A\|^d$.

Inequality~\eqref{e.1st ineq} is a very simple consequence of the Cayley-Hamilton theorem.
Indeed, let 
$p(z) = z^d - \sigma_1 z^{d-1} + \cdots + (-1)^d \sigma_d$ be the
characteristic polynomial of $A$.
Since $p(A) = 0$, 
$$
\|A^d\| \le \sum_{i=1}^d |\sigma_i| \; \|A\|^{d-i}.
$$
Since the $\sigma_i$ are the elementary symmetric functions on the eigenvalues of~$A$,
$$
|\sigma_i| \le \binom{d}{i} \rho(A)^i \le \binom{d}{i} \rho(A) \; \|A\|^{i-1}.
$$
Therefore~\eqref{e.1st ineq} follows.

\smallskip

The \emph{spectral radius theorem} (for the finite-dimensional case) asserts that
\begin{equation}\label{e.srf}
\rho(A) = \lim_{n \to \infty} \| A^n \|^{1/n}.
\end{equation}
The formula above may be deduced from the inequality~\eqref{e.1st ineq}, 
as we now show.
Since $\|A^{n+m}\| \le \|A^n\| \, \|A^m\|$,
the limit in~\eqref{e.srf} exists (see \mbox{\cite[problem I.98]{Polya Szego}});
let us call it $r$.
Clearly, $r \ge \rho(A)$.
Applying~\eqref{e.1st ineq} to $A^n$ in the place of $A$,
using that $\rho(A^n) = \rho(A)^n$ and taking the $1/dn$-power, we obtain
$$
\| A^{dn} \|^{1/dn} \le C^{1/dn} \rho(A)^{1/d} \| A^n \|^{(d-1)/dn}.
$$
Taking limits when $n \to \infty$, we get
$r \le \rho(A)^{1/d} \; r^{(d-1)/d}$,
that is, $r \le \rho(A)$, proving~\eqref{e.srf}.
The author ignores whether this proof
has ever appeared in the literature.

\medskip

Now, let $\Sigma$ be a non-empty bounded subset of $\md$.
Define
$$
\| \Sigma \| = \sup_{A \in \Sigma} \|A\|, \qquad
\rho(\Sigma) = \sup_{A \in \Sigma} \rho(A).
$$
If $n \in \N$, we denote by $\Sigma^n$ the
set of the products $A_1 \cdots A_n$, with all $A_i \in \Sigma$.
Since $\|\Sigma^{n+m}\| \le \|\Sigma^n\| \, \|\Sigma^m\|$, the limit 
$$
\RR(\Sigma)= \lim_{n \to \infty} \| \Sigma^n \|^{1/n}
$$
exists and equals $\inf_n \| \Sigma^n \|^{1/n}$.
Besides, it is independent of the chosen norm.
The quantity $\RR(\Sigma)$ was introduced by Rota and Strang \cite{Rota Strang}
and is called the \emph{joint spectral radius} of the set $\Sigma$.
For a nice geometrical interpretation of the joint spectral radius,
see~\cite{Protasov Russian} (or~\cite{Protasov}).

\smallskip

Our first main result is a generalization of~\eqref{e.srf} to sets of matrices:

\begin{mr}\label{t.a}
Given $d \ge 1$ , there exists $C_1>1$ such that,
for every bounded set $\Sigma \subset \md$ and every norm $\norma$ in $\C^d$,
$$
\| \Sigma^d \| \le C_1  \RR(\Sigma)  \| \Sigma \|^{d-1}.
$$
\end{mr}

Our next result relates the joint spectral radius of $\Sigma$
with spectral radii of products of matrices in $\Sigma$:

\begin{mr}\label{t.b}
Given $d \ge 1$, there exists $C_2>1$ and $k \in \N$ such that
for every bounded set $\Sigma \subset \md$,
$$
\RR(\Sigma) \le C_2 \max_{1 \le j \le k} \rho(\Sigma^j)^{1/j}.
$$
\end{mr}

Using theorem~\ref{t.b}, we can extend the spectral radius theorem~\eqref{e.srf}:

\begin{corollary}[Berger-Wang generalized spectral radius theorem] \label{c.gsrt}    
If $\Sigma \subset \md$ is bounded then 
$$
\RR(\Sigma) = \limsup_{n \to \infty} \rho (\Sigma^n) ^{1/n}.
$$
\end{corollary}

\begin{proof}
The inequality $\RR(\Sigma) \ge \limsup \rho (\Sigma^n) ^{1/n}$ is trivial.
Applying theorem~\ref{t.b} to $\Sigma^n$ and using that $\RR(\Sigma^n) = \RR(\Sigma)^n$,
we obtain
$$
\RR(\Sigma) \le C_2^{1/n} \max_{1 \le j \le k} \rho(\Sigma^{jn})^{1/{jn}}.
$$
Taking $\limsup$ when $n\to \infty$, we get the result.
\end{proof}

The result above was conjectured by Daubechies and Lagarias~\cite{Daubechies Lagarias}
and proved by Berger and Wang~\cite{Berger Wang}.
Other proofs were given in~\cite{Elsner} and \cite{Shih Wu Pang}.

\smallskip

The proof of theorem~\ref{t.a} is elementary, while
in the proof of theorem~\ref{t.b} we shall use some geometric invariant theory
We also give another generalization of~\eqref{e.1st ineq}, proposition~\ref{p.other ineq} below,
whose proof is elementary.

%\begin{remark}
%We shall actually prove an inequality that is stronger than theorem~\ref{t.b}:
%$\RR(\Sigma) \le C_2
%\sup \{ \, \lvert \tr(A_1 \cdots A_j) \rvert^{1/j} ; \, 
%1 \le j \le k, \, A_1,\ldots,A_j \in \Sigma \, \}$.
%\end{remark}

\begin{remark}
For all $\Sigma$ and $m$, $n\in\N$, we have
$\rho(\Sigma^{mn})^{1/mn} \ge \rho(\Sigma^n)^{1/n}$
(because $\Sigma^{nm} \subset (\Sigma^n)^m$).
So in theorem~\ref{t.b} it is sufficient to take the maximum of $\rho(\Sigma^j)^{1/j}$
over $j$ with $k/2 < j \le k$.
Another consequence of the latter remark is that
$$
\limsup_{n \to \infty} \rho (\Sigma^n) ^{1/n} = \sup_{n \in \N} \rho (\Sigma^n) ^{1/n}.
$$
\end{remark}

%%%%%%%%%%%%%%%%%%%%%%%%%%%%%%%%%%%%%%%%%%%%%%%%%%%%%%%%

\section{Proof of theorem~\ref{t.a}}

We first prove an inequality that is weaker than theorem~\ref{t.a}:

\begin{lemma}\label{l.weaker ineq}
Let $\norma_e$ be the euclidian norm in $\C^d$.
There exists $C_0 = C_0(d)$ such that 
$$
\| S \Sigma^d S^{-1} \|_e \le C_0 \| \Sigma \|_e \, \| S \Sigma S^{-1} \|_e^{d-1} .
$$
for every non-empty bounded set $\Sigma \subset \md$ and every $S \in \gldc$.
\end{lemma}

\begin{proof}
We shall also consider the norm in $\md$ defined by
$$
\|A\|_0 = \max |a_{ij}|,
\text{ where } A = (a_{ij})_{i,j=1,\ldots,d}.
$$
We first assume $S$ is a diagonal matrix $\mathrm{diag}(\lambda_1, \ldots, \lambda_d)$,
with $\lambda_1, \ldots, \lambda_d > 0$.
Take $d$ matrices $A_1, \ldots, A_d \in \Sigma$, and write $A_\ell = (a^{(\ell)}_{ij})$.
Then
$$
\| S \Sigma S^{-1} \|_0 =
\max_{i,j,\ell} |\lambda_i a^{(\ell)}_{i j} \lambda_j^{-1}|.
$$
and
$$
\| S A_1 \cdots A_d S^{-1} \|_0 \le
C_0 \max_{i_0,\ldots, i_d} 
\big| \lambda_{i_0} a^{(1)}_{i_0 i_1} \cdots a^{(d)}_{i_{d-1} i_d} \lambda_{i_d}^{-1} \big|,
$$
where $C_0 = d^{d-1}$.
Given integers $i_0, \ldots, i_d \in \{1,\ldots,d \}$, 
by the pigeon-hole principle there exists
$1 \le k \le d$ such that $\lambda_{i_{k-1}} \le \lambda_{i_k}$.
Therefore
\begin{multline*}
\big| \lambda_{i_0} a^{(1)}_{i_0 i_1} \cdots a^{(d)}_{i_{d-1} i_d} \lambda_{i_d}^{-1} \big| =
\prod_{1 \le \ell \le d} | \lambda_{i_{\ell-1}} a^{(\ell)}_{i_{\ell-1} i_\ell} \lambda_{i_\ell}^{-1} | 
\le \\ \le
|a^{(k)}_{i_{k-1},i_k}|
\prod_{\ell \neq k} | \lambda_{{i_{\ell-1}}} a^{(\ell)}_{i_{\ell-1} i_\ell} \lambda_{i_\ell}^{-1} |
\le
\|\Sigma \|_0 \, \| S \Sigma S^{-1} \|_0^{d-1} .
\end{multline*}
It follows that 
$\| S \Sigma^d S^{-1} \|_0 \le
C_0 \| \Sigma \|_0 \, \| S \Sigma S^{-1} \|_0^{d-1}$.
Up to changing $C_0$, the same inequality holds for the euclidian norm $\norma_e$.

Next consider the general case $S \in \gldc$.
By the singular value decomposition theorem,
there exist unitary matrices $U$, $V$ and a diagonal matrix 
$D = \mathrm{diag}(\lambda_1, \ldots, \lambda_d)$, with
$\lambda_1, \ldots, \lambda_d > 0$, such that $S = UDV$. 
Since $U$ and $V$ preserve the euclidian norm,
\begin{multline*}
\| S \Sigma^d S^{-1} \|_e 
  =   \| D (V \Sigma V^{-1})^d D^{-1} \|_e \le \\
  \le C_0 \| V \Sigma V^{-1}\|_e \, \| DV \Sigma V^{-1} D^{-1} \|_e^{d-1} 
  =   C_0 \| \Sigma \|_e \, \| S \Sigma S^{-1} \|_e^{d-1}.
\end{multline*}
This proves the lemma.
\end{proof}

To make the constant in lemma~\ref{l.weaker ineq} independent of the norm, we will~use:

\begin{lemma} \label{l.equiv norms}
There exists $C = C(d)$ such that, for every two norms $\norma_1$, $\norma_2$ in $\C^d$,
there is $S \in \gldc$ such that:
\begin{enumerate}
\item $C^{-1} \|v\|_1 \le \|Sv\|_2\le \|v\|_1$, for all $v \in \C^d$;
\item $C^{-1} \|A\|_1 \le \|SAS^{-1}\|_2 \le C \|A\|_1$, for all $A \in \md$.
\end{enumerate}
\end{lemma}

\begin{proof}
The second part is an immediate consequence of the first one.
To prove the first part, it is enough to show that
for every $\norma$ in $\C^d$, there is $S \in \gldc$
such that 
\begin{equation}\label{e.claim}
C_d^{-1} \|v\| \le \| Sv \|_0 \le \|v\|, \quad \forall v \in \C^d,
\end{equation}
where $\norma_0$ is the sup-norm in $\C^d$ and $C_d = 2^d -1$.
The proof is by induction. 
Let $\norma$ be a norm in $\C^{d+1}$.
Restrict it to the subspace $\C^d  = \C^d \times\{0\}\subset \C^{d+1}$.
By induction hypothesis, there is $S\in \gldc$ such that~\eqref{e.claim} holds.
%%$S:\C^d \to \C^{d+1}$ such that $S(\C^d) = \C^d$ and~\eqref{e.claim} holds.
If $\pi_j : \C^{d+1} \to \C$ is the projection in the $j$-th coordinate,
then $\lvert \pi_j \circ S(v) \rvert \le \|v\|$, for all $v \in \C^d$ and $1 \le j \le d$.
By the Hahn-Banach theorem, there are linear functionals
$\Lambda_j: \C^{d+1} \to \C$ such that $\Lambda_j|\C^d = \pi_j \circ S$ and 
$\lvert \Lambda_j(w) \rvert \le \|w\|$, for all $w \in \C^{d+1}$ and $1 \le j \le d$.
Let $a= \| \pi_{d+1} \|$ and
define a linear map $\bar{S}: \C^{d+1} \to \C^{d+1}$ by
$$
\pi_j\circ\bar{S} =  
\begin{cases}
\Lambda_j        &\text{if }1 \le j \le d,\\
a^{-1} \pi_{d+1} &\text{if }j=d.
\end{cases}
$$
Then $\bar{S}\in GL(d+1)$ and $\|\bar{S}w\|_0 \le \|w\|$,
so $\bar{S}$ satisfies the second inequality in~\eqref{e.claim}.
To prove the first one, 
let $\xi \in \C^{d+1}$ be such that $\pi_{d+1}(\xi)=a$ and $\|\xi\|=1$.
Write $\bar{S}(\xi) = \eta + e_{d+1}$ with $\eta \in \C^d$ and 
$e_{d+1}=(0,\ldots,0,1) \in \C^{d+1}$.
We have $\|\eta\|_0 \le \|\xi\| = 1$ and so
$\| \bar{S}^{-1}(\eta) \| = \| S^{-1}(\eta) \| \le C_d \|\eta\|_0 \le C_d$.
Therefore
$$
\| \bar{S}^{-1}(e_{d+1}) \| \le \|\xi\| + \| \bar{S}^{-1}(\eta) \| \le 1+ C_d.
$$
Now let $w \in \C^{d+1}$ be given.
Write $w = v + t e_{d+1}$ with $w\in\C^d$ and $t\in\C$.
Then
\begin{multline*}
\| \bar{S}^{-1} (w) \| 
 \le \| \bar{S}^{-1} (v) \| + |t| \, \| \bar{S}^{-1} (e_{d+1}) \|
 \le C_d \|v\|_0 + (C_d + 1)|t|  \\
 \le (2C_d+1) \max\{\|v\|_0,|t|\}
 = C_{d+1} \|w\|_0. 
\end{multline*}
This proves that~\eqref{e.claim} holds with $d+1$ and $\bar{S}$
in the place of $d$ and $S$.
\end{proof}

The result below gives another characterization of the joint spectral radius.
For a proof, see~\cite{Elsner} or \cite{Rota Strang}.

\begin{proposition} \label{p.rs}
For all bounded $\Sigma \subset \md$,
$$
\RR(\Sigma) = \inf_{\norma} \|\Sigma\|,
$$
where the infimum is taken over all norms in $\C^d$.
\end{proposition}

\begin{proof}[Proof of theorem~\ref{t.a}]
Let $C_0$ and $C$ be as in lemmas~\ref{l.weaker ineq} and \ref{l.equiv norms}.
Let $\norma_e$ be the euclidian norm, and let $\norma_1$, $\norma_2$ 
be any two norms in $\C^d$.
Let $S_1$, $S_2 \in \gldc$ be given by lemma~\ref{l.equiv norms} such that
$$
C^{-1} \| S_i A S_i^{-1}\|_e \le \|A\|_i \le C \| S_i A S_i^{-1}\|_e \quad
\forall A \in \md, \, i=1,2.
$$
Take $\Sigma \subset \md$.
Then (applying lemma~\ref{l.weaker ineq} with $S = S_1 S_2^{-1}$ and 
$S_2 \Sigma S_2^{-1}$ in the place of $\Sigma$) 
$$
\|\Sigma^d\|_1 
   \le C       \| S_1 \Sigma^d S_1^{-1} \|_e 
   \le C C_0   \| S_2 \Sigma S_2^{-1} \|_e \| S_1 \Sigma S_1^{-1} \|_e^{d-1}
   \le C^d C_0 \| \Sigma \|_2 \| \Sigma \|_1^{d-1}.
$$
Taking the infimum over $\norma_2$ in the left hand side,
we obtain, by proposition~\ref{p.rs},
$\|\Sigma^d\|_1  \le C_1 \RR(\Sigma) \|\Sigma\|_1^{d-1}$,
where $C_1 = C^d C_0$.
\end{proof}

Let us reread theorem~\ref{t.a} in terms of another invariant.
Given a non-empty bounded $\Sigma \subset \md$, we define
\begin{equation} \label{e.def SS}
\SS(\Sigma) = \sup_{\norma} \frac{\|\Sigma^d\|}{\|\Sigma\|^{d-1}}
\qquad \text{if $\Sigma \neq \{0\}$,}
\end{equation}
and $\SS(\{0\}) = 0$.
The functions $\RR(\mathord{\cdot})$ and $\SS(\mathord{\cdot})$ are comparable:

\begin{proposition} \label{p.RR SS}
$\RR(\Sigma) \le \SS(\Sigma) \le C_1 \RR(\Sigma)$.
\end{proposition}

\begin{proof}
The second inequality is theorem~\ref{t.a}.
For any $\norma$ we have $\|\Sigma^d\| \ge \RR(\Sigma)^d$ and
so, using proposition~\ref{p.rs}, 
$$
\SS(\Sigma) \ge \sup_{\norma} \frac{\RR(\Sigma)^d}{\| \Sigma \|^{d-1}} = 
\RR(\Sigma).
$$
\end{proof}

%%%%%%%%%%%%%%%%%%%%%%%%%%%%%%%%%%%%%%%%%%%%%%%%%%%%%%%%
\section{Proof of theorem~\ref{t.b}}

We shall need the following general result:

\begin{proposition} \label{p.invariant bound}
Fix $d$, $\ell \ge 1$.
Let $f: \md^\ell \to [0, \infty)$ be a locally bounded function such that,
for every $A_1, \ldots, A_\ell \in \md$,
\begin{itemize}
\item
$f (S A_1 S^{-1}, \ldots, S A_\ell S^{-1}) = f (A_1, \ldots, A_\ell)$,
$\forall S \in \gldc$;
\item
$f (t A_1, \ldots, t A_\ell) = |t| f ( A_1, \ldots,  A_\ell)$, 
$\forall t \in \C$.
\end{itemize}
Then there exist numbers $k=k(d) \in \N$ and $C=C(d,\ell, f)>0$ such that
\begin{equation} \label{e.result}
f (A_1,\ldots,A_\ell) \le C \max_{1 \le j \le k} \rho(\Sigma^j)^{1/j},
\qquad \text{where } \Sigma = \{ A_1, \ldots, A_\ell \}.
\end{equation}
\end{proposition}

Let us postpone the proof of this proposition and conclude the:

\begin{proof}[Proof of theorem~\ref{t.b}]
Let $\SS(\mathord{\cdot})$ be as in~\eqref{e.def SS}.
Define a function $f: \md^d \to [0,\infty)$ by
$f(A_1, \ldots, A_d) = \SS(\{A_1, \ldots, A_d\})$.
By theorem~\ref{t.a}, $f(\Sigma) \le C_1 \| \Sigma \|$ (for any norm) --
in particular, $f$ is locally bounded.
$f$ also satisfies the other hypotheses of proposition~\ref{p.invariant bound},
thus there are $k$ and $C_2$ such that
\begin{equation}\label{e.SS bound}
\SS (\Sigma) \le C_2 \max_{1 \le j \le k} \rho(\Sigma^j)^{1/j},
\end{equation}
for every $\Sigma \subset \md$ with at most $d$ elements.
But 
$$
\SS(\Sigma) = \sup\{\SS(\Sigma'); \, \Sigma' \subset \Sigma, \,  \# \Sigma' \le d  \},
$$
hence~\eqref{e.SS bound} actually holds for every bounded $\Sigma$.
Since $\RR(\Sigma) \le \SS(\Sigma)$ (proposition~\ref{p.RR SS}),
theorem~\ref{t.b} follows.
\end{proof}

A few preliminaries in geometric invariant theory 
are necessary to prove proposition~\ref{p.invariant bound}.
Some references are~\cite{Newstead} and~\cite{Kraft}.

\subsection{Polynomial invariants}
Let $V$ be a complex vector space,
$G$ be a group and $\iota: G \to GL(V)$ be a linear representation of $G$.
We shall write $gx = \iota(g)(x)$.
The \emph{orbit} of $x\in V$ is the set $\OO(x) = \{ gx; \; g \in G \}$. 
Let $\C[V]$ be the ring of polynomial functions $\phi: V \to \C$.
A polynomial $\phi \in \C[V]$ is \emph{invariant} if it is constant along
each orbit, that is, $\phi(g x) \equiv \phi(x)$.
The \emph{ring of invariants}, denoted by $\C[V]^G$,
is the set of all invariant polynomials.

For some groups $G$, called \emph{reductive} groups,
a celebrated theorem of Nagata asserts that the ring $\C[V]^G$ is finitely generated.
We shall not define a reductive group; but some examples are $\gldc$, $\sldc$, $\pgldc$.
We assume from now on that $G$ is reductive.
In this case, the theory provides an algebraic quotient of $V$ by $G$
with good properties:

\begin{theorem}\label{t.quotient}
Let $\phi_1, \ldots, \phi_N$ be a set of generators of $\C[V]^G$.
Let $\pi: V \to \C^N$ be the mapping
$x \mapsto (\phi_1(x), \ldots, \phi_N(x)) \in \C^N$.
Then:
\begin{itemize}
\item[0.] $\pi$ is $G$-invariant (i.e., constant along orbits);
\item[1.] $Y = \pi(V)$ is closed;
\item[2.] $\pi(x_1) = \pi(x_2)$ if and only if the closures 
$\overline{\OO(x_1)}$ and $\overline{\OO(x_2)}$ have non-empty intersection; 
\item[3.] for every $y\in Y$, the fiber $\pi^{-1}(y)$ contains an
unique closed orbit.
\end{itemize}
\end{theorem}

In the statement above, and in everything that follows, 
the spaces $V$ and $\C^N$ are endowed with the ordinary (not Zariski) topologies.
Notice item 2 says that $\pi$ separates every pair of orbits that 
can be separated by a $G$-invariant continuous function.

\begin{proof}[Indication of proof]
Let $\C(V)^G$ be the field of $G$-invariant rational functions.
It's easy to see that $\C(V)^G$ is the field of quotients of $\C[V]^G$.
Let $\pi$ and $Y$ be as in the statement.
Let $Z$ be the Zariski-closure of $Y$, and consider $\pi$ as a function $V \to Z$.
Then $\pi$ induces a homomorphism $\pi^*:C(Z) \to \C(V)$ via $f \mapsto f \circ \pi$.
One easily shows that $\pi^*$ is an isomorphism onto $\C(V)^G$, so
$\pi:V \to Z$ is an algebraic quotient in the sense of~\cite[\S II.3.2]{Kraft}.
Therefore $\pi$ is surjective, that is, $Y = Z$, and item 1 follows.
Items 2 and 3 are \cite[corollary~3.5.2]{Newstead} and \cite[bemerkung~1,~\S II.3.2]{Kraft},
respectively.
In the references above the Zariski topology is used instead.
But this makes no difference here, by \cite[\S AI.7]{Kraft}.
\end{proof}

\begin{example}
Let $G = \gldc$ act on $V = \md$ by conjugation: 
$\iota(S)(A) = SAS^{-1}$ for $S \in G$ and $A\in V$.
Given $A\in V$, let $\sigma_1(A), \ldots, \sigma_d(A)$ be the coefficients of
the characteristic polynomial of $A$.
Then $\sigma_1, \ldots, \sigma_d \in \C[V]^G$.
Moreover, these polynomials generate the ring $\C[V]^G$.
Let $\pi = (\sigma_1, \ldots, \sigma_d)$.
Then $\pi$ is onto $\C^d$.
Every fiber $\pi^{-1}(y)$ consists in finitely many orbits,
each one corresponding to a different Jordan form.
The closed orbits are those of diagonalizable matrices.
The fiber $\pi^{-1}(0)$ is the set of nilpotent matrices.
(See \cite[\S I.3]{Kraft}.)
\end{example}

\subsection{Topological considerations}
Fix $\pi$ and $Y = \pi(V)$ as in theorem~\ref{t.quotient}.
We endow the set $Y \subset \C^N$ with the induced topology.
The following theorem was proved independently by 
Luna~\cite{Luna} and Neeman~\cite[corollary~1.6, remark~1.7]{Neeman}
(see also \cite{NeemanSurvey}):

\begin{theorem} \label{t.topology}
The topology in $Y$ coincides with the quotient topology induced by $\pi: V \to Y$
(i.e., $U \subset Y$ is open if and only if $\pi^{-1}(U)$ is open in $V$).
\end{theorem}

\begin{corollary} \label{c.semiproper}
The mapping $\pi: V \to Y$ is \emph{semiproper}, that is,
for every compact set $L \subset Y$ there exists a
compact set $K \subset V$ such that $\pi(K) \supset L$.
\end{corollary}

\begin{proof}
Suppose that for some compact $L \subset Y$ 
there is no compact set \mbox{$K \subset V$} such that $\pi(K) \supset L$.
Take compact sets $K_n \subset V$ such that
$K_n \subset \mathrm{int}\, K_{n+1}$ and $\bigcup_n K_n = V$.
Then, for each $n$, there exists $y_n \in L$ such that 
$\pi^{-1}(y_n) \cap K_n = \varnothing$.
Up to replacing $(y_n)$ with a subsequence, we may assume
that $y = \lim y_n$ exists and $y_n \neq y$ for each $n$.
Then the set $F = \{ y_n ; \; n \in \N \}$ is not closed in $Y$,
but $\pi^{-1}(F) = \bigcup_n \pi^{-1}(y_n)$ is closed in $V$,
contradicting theorem~\ref{t.topology}.
\end{proof}

Let us derive a consequence of the above results:

\begin{lemma} \label{l.goes down}
If $f:V \to [0,\infty)$ is a $G$-invariant locally bounded function then
there exists a locally bounded $h:Y \to [0,\infty)$ such that $f \le h \circ \pi$.
\end{lemma}

\begin{proof}
Given $x \in V$, let $F_x = \pi^{-1}(\pi(x))$ be the fiber containing $x$.
Set
$$
\bar{f}(x) = \inf 
\big\{ \sup f|U ; \; \text{$U$ is a $G$-invariant open set containing $F_x$ } \big\}.
$$
(Here ``$U$ is $G$-invariant'' means $\OO(x) \subset U$ for all $x \in U$.)
We claim that $\bar{f}(x)$ is finite for all $x \in V$.
Indeed, each fiber $F_x$ contains an unique closed orbit $\OO(x_0)$, by theorem~\ref{t.quotient}.
Let $U_0$ be a bounded neighborhood of $x_0$; so $\sup f|U_0$ is finite.
Let $U = \bigcup_{x \in U_0} \OO(x)$;
then $U$ is a $G$-invariant open set and $\sup f|U = \sup f|U_0$.
Moreover, $U$ contains $F_x$:
for every $\xi \in F_x$, we have, by theorem~\ref{t.quotient},
$\overline{\OO(\xi)} \cap \OO(x_0) \neq \varnothing$,
hence $\OO(\xi) \cap U_0 \neq \varnothing$ and $\xi\in U$.
This proves that $\bar{f}(x) \le \sup f|U < \infty$.

The function $\bar{f}: V \to \R$ satisfies $\bar{f} \ge f$ and is also locally bounded.
Since $\bar{f}$ is constant on fibers, there exist $h:Y \to \R$ such that
$\bar{f} = h \circ \pi$.
The function $h$ is locally bounded, because if $L \subset Y$ is a compact set
then, by corollary~\ref{c.semiproper}, there is some compact $K \subset V$
such that $\pi(K) \supset L$ and, in particular, 
$h|L \le (h \circ \pi)|K = \bar{f}|K < \infty$.
\end{proof}

\subsection{$\ell$-uples of matrices and end of the proof}

From now on we set $G=\gldc$, $V=\md^d$ and
$$
\iota(S) (A_1, \ldots, A_\ell) = (S A_1 S^{-1}, \ldots, S A_\ell S^{-1}).
$$
In this case, a finite set of generators for $\C[V]^G$ is known:

\begin{theorem}[Procesi \cite{Procesi}, theorem 3.4a] \label{t.procesi}
The ring of invariants is generated by the polynomials
$\tr (A_{i_1} \cdots A_{i_j})$, with $1 \le j \le k$,
where $k=2^d -1$.
\end{theorem}

We are now able to give the:
\begin{proof}[Proof of proposition~\ref{p.invariant bound}]
Let $k = 2^d -1$, $N = \ell + \ell^2 + \cdots + \ell^k$, and let
$\alpha_1, \ldots, \alpha_N$ be all the sequences
$\alpha = (i_1, \ldots, i_j) \in \{1,\ldots,\ell\}^j$ of length $|\alpha| =j$,
$1 \le j \le k$.
Let $\pi = (\phi_1, \ldots, \phi_N): V \to \C^N$ be given by 
$$
\phi_i (A_1,\ldots,A_\ell) = \tr (A_{i_1} \cdots A_{i_j}), \quad
\text{where } \alpha_i = (i_1, \ldots, i_j).
$$
Let $Y = \pi(V)$.
Define another function $\tau: \C^N \to \R$ by
$$
\tau(z_1, \ldots, z_N) = 
\max \{ |z_\ell|^{1/|\alpha_i|} ; \; 1 \le i \le N \}.
$$
So if $x = (A_1, \ldots, A_\ell)$ then
$$
\tau(\pi(x)) = \max \big\{ \lvert \tr (A_{i_1} \cdots A_{i_j}) \rvert^{1/j} ; \;
1 \le j \le k, \; 1 \le i_1, \ldots, i_j \le \ell \big\}.
$$
Since $\lvert \tr(A) \rvert \le d \rho(A)$ for every $A\in \md$, we have
$$
\tau(\pi(x)) \le d \max_{1\le j \le k} \rho(\Sigma^j)^{1/j},
\qquad \text{where } \Sigma = \{ A_1, \ldots, A_\ell \}.
$$
Notice $\tau(\pi(tx)) = |t| \tau (\pi(x))$, for all $t\in\C$.
Let $h:Y \to [0,\infty)$ be given by lemma~\ref{l.goes down}.
Since $K= \tau^{-1}(1)$ is compact and $Y$ is closed, $C_0 = \sup h|(Y \cap K)$ is finite.
Given $x \in V$, let $t = \tau(\pi(x))$.
If $t\neq 0$ then
$$
f(x) = t f(t^{-1} x) \le t h (\pi(t^{-1} x)) \le C_0 t = C_0 \tau(\pi(x))
\le d C_0 \max_{1\le j \le k} \rho(\Sigma^j)^{1/j}.
$$
Let $C = d C_0$; then~\eqref{e.result} holds.
If $t=0$, that is, $\pi(x)=0$, we argue differently.
By theorem~\ref{t.quotient}, the orbit of $x$ accumulates at $0$.
It follows from the hypotheses on $f$ that $f(0)=0$ and $f$ is continuous at $0$.
Therefore $f(x)=0$ and~\eqref{e.result} holds.
This completes the proof of proposition~\ref{p.invariant bound}
and so of theorem~\ref{t.b}.
\end{proof}

%%%%%%%%%%%%%%%%%%%%%%%%%%%%%%%%%%%%%%%%%%%%%%%%%%%%%%%%
\section{Another inequality and some questions}

We shall prove another inequality, proposition \ref{p.other ineq} below,
which generalizes~\eqref{e.1st ineq} and is also an elementary consequence of the Cayley-Hamilton theorem.

We need some notation.
For $s\ge 1$, let $S_s$ be the set of permutations of $\{1,2,\ldots,s\}$.
Given $\sigma \in S_s$, decompose $\sigma$ in disjoint cycles, including the ones of length $1$:
$$
\sigma = (i_1 \cdots i_k) (j_1 \cdots j_h) \cdots (t_1 \cdots t_e).
$$
Then, given matrices $A_1,\ldots,A_s$, we set 
$$
\Phi_\sigma(A_1,\ldots,A_s) =
\tr(A_{i_1} \cdots A_{i_k}) \tr(A_{j_1} \cdots A_{j_h}) \cdots \tr(A_{t_1} \cdots A_{t_e}).
$$
Letting $\eps(\sigma)$ be the sign of $\sigma$, we define
\begin{equation}\label{e.def F}
F(A_1,\ldots,A_s) = \sum_{\sigma \in S_s} \eps(\sigma) \Phi_\sigma(A_1,\ldots,A_s).
\end{equation}
Define also
$P(A_1,\ldots,A_s)=
\sum_{\sigma \in S_s} A_{\sigma(1)} \cdots A_{\sigma(s)}$.
The \emph{trace identity} from \cite[Corollary 4.4]{Procesi}
(which follows from the Cayley-Hamilton theorem by an elementary process, 
see also \cite[\S 4]{Formanek}) is:
\begin{equation}\label{e.trace id}
\sum_{s=0}^d \sum (-1)^s F(A_{i_1}, \ldots, A_{i_s}) P(A_{j_1}, \ldots, A_{j_{d-s}})  = 0 ,
\end{equation}
where the second sum runs over all partitions of $\{1,\ldots,d\}$ into
two disjoint subsets $\{i_1 < \ldots < i_s\}$ and $\{j_1 < \ldots < j_{d-s}\}$;
it is understood that $F(\varnothing) = 1$ and  $P(\varnothing) = I$.

\begin{proposition}\label{p.other ineq}
Given $d \ge 1$, there exists $C>1$ such that 
for every operator norm $\norma$ and every $d$ matrices
$A_1,\ldots,A_d \in \md$, we have
$$
\| P(A_1, \ldots A_d) \| \le
C \| \Sigma \|^{d-1} \max_{1 \le j \le d} \rho(\Sigma^j)^{1/j}.
$$
where $\Sigma = \{ A_1, \ldots, A_d \}$.
\end{proposition}

\begin{proof}
We estimate terms in \eqref{e.trace id}, for $1 \le s \le d$.
If $\sigma$ is a permutation of $\{i_1 < \ldots < i_s\}$ with cycles of lengths $k_1, \ldots, k_h$ then
$$
|\Phi_\sigma (A_{i_1}, \ldots, A_{i_s})| \le C_0 \rho(\Sigma^{k_1}) \cdots \rho(\Sigma^{k_h}),
$$
where $C_0$ is a constant.
The right hand side is $\le C_0 \rho(\Sigma^{k_i})^{1/k_i} \|\Sigma\|^{s-1}$, for any $k_i$.
Plugging this estimate in \eqref{e.def F}, we get
$$
|F(A_{i_1}, \ldots, A_{i_s})| \le C_0 \|\Sigma\|^{s-1} \max_{1 \le j \le s} \rho(\Sigma^j)^{1/j}. 
$$
Using the inequality above and 
the obvious bound $\|P(A_{j_1}, \ldots, A_{j_{d-s}})\| \le (d-s)! \|\Sigma\|^{d-s}$,
the result follows from \eqref{e.trace id}.
\end{proof}

We do not know whether the methods of the proof of proposition~\ref{p.other ineq}
can be improved to give an elementary proof of theorem~\ref{t.b}.
Notice that if $k$ in theorem~\ref{t.b} were equal to $d$ then 
proposition~\ref{p.other ineq} would follow from theorems~\ref{t.a} and \ref{t.b}.

\begin{question}
What is the minimum $k$ such that theorem~\ref{t.b} holds? %\margem{Nagata-Higman?}
Can one take $k=d$?
\end{question}

The answer is yes when $d=2$.
The ring of invariants of two $2 \times 2$ matrices $A_1$ and $A_2$ is generated by
$\tr A_1$, $\det A_1$, $\tr A_2$, $\det A_2$, $\tr A_1 A_2$,
see \cite[\S 7]{Formanek}.
Since $\det A$ can be expressed as a polynomial in $\tr A$ and $\tr A^2$,
one can take $k=2$ in theorem~\ref{t.procesi}, and so also in theorem~\ref{t.b},
when $d=2$.
Moreover, since $\rho(\Sigma) \le \rho(\Sigma^2)^{1/2}$,
theorem~\ref{t.b} assumes the form:
$$
\RR(\Sigma) \le C_2 \rho(\Sigma^2)^{1/2}.
$$
Using this inequality, it is easy to show that the sequence $\rho(\Sigma^{2n})^{1/2n}$ converges.
However, the sequence $\rho(\Sigma^n)^{1/n}$ itself does not necessarily converge.
We reproduce an example from \cite{Gripenberg}:
$$
\Sigma = \left\{
\begin{pmatrix}
0 & 1 \\ 0 & 0
\end{pmatrix},
\begin{pmatrix}
0 & 0 \\ 1 & 0
\end{pmatrix}
\right \}
\; \Longrightarrow \; 
\rho(\Sigma^n) =
\begin{cases}
0 &\text{if $n$ is odd,}\\
1 &\text{if $n$ is even.}
\end{cases}
$$

%%%%%%%%%%%%%%%%%%%%%%%%%%%%%%%%

\begin{ack}
I would like to thank J.~V.~Pereira for helpful discussions.
\end{ack}

%%%%%%%%%%%%%%%%%%%%%%%%%%%%%%%%

\end{document}